\documentclass[12pt, a4paper]{article}
\usepackage[latin1]{inputenc}
\usepackage[english]{babel}
\usepackage{amsmath}
\usepackage[english]{babel}
\usepackage{indentfirst}
\usepackage{amstext}
\usepackage{enumerate}
\usepackage{amsfonts}
\usepackage{textcomp}
\usepackage{amssymb}
\usepackage[dvips]{graphicx}
\usepackage{setspace}
\usepackage[all]{xy}
\usepackage{color}

\newcommand {\demo}{\hskip -0.6cm{\bf Proof.  }}
\newcommand {\fim}{\hfill{$\square$}\vskip 1pc}

\newcommand {\N}{\mathbb{N}}
\newcommand {\Z}{\mathbb{Z}}

\newcommand {\F}{\mathbb{F}}

\newcommand {\GG}{\mathcal{G}}

\newtheorem{teorema}{Theorem}[section]
\newtheorem{lema}[teorema]{Lemma}
\newtheorem{corolario}[teorema]{Corollary}
\newtheorem{definicao}[teorema]{Definition}
\newtheorem{proposicao}[teorema]{Proposition}

\newtheorem{remark}[teorema]{Remark}
\newtheorem{notation}[teorema]{Notation}

\begin{document}
\onehalfspace

\title{Simplicity and chain conditions for ultragraph Leavitt path algebras via partial skew group ring theory}
\maketitle
\begin{center}

{\large Daniel Gon\c{c}alves\footnote{This author is partially supported by CNPq.} and Danilo Royer}\\
\end{center}  
\vspace{8mm}

\begin{abstract}
We realize Leavitt ultragraph path algebras as partial skew group rings. Using this realization we characterize artinian ultragraph path algebras and give simplicity criteria for these algebras. 
\end{abstract}

\vspace{1.5pc}

\vspace{0.5pc}


\vspace{1.5pc}
\section{Introduction}

The study of algebras associated to combinatorial objects has attracted a great deal of attention in the past years. Part of the interest in these algebras arise from the fact that many properties of the combinatorial object translate into algebraic properties of the associated algebras and, furthermore, there are deep connections between these algebras and symbolic dynamics. As examples of algebras associated to combinatorial objects we cite graph C*-algebras, Leavitt path algebras, higher rank graph algebras, Kumjian-Pask algebras, ultragraph C*-algebras, among others (see \cite{GA, AASbook} for a comprehensive list). 

Notice that in the list of algebras we presented above the C*-algebraic version of the algebras was immediately followed by the algebraic analogue, except for the ultragraph case. Ultragraphs (a generalization of graphs, where the range map takes values on the power set of the vertices) were defined by Mark Tomforde in \cite{Tom3} as an unifying approach to Exel-Laca and graph C*-algebras. They have proved to be a key ingredient in the study of Morita equivalence of Exel-Laca and graph C*-algebras (see \cite{KMST}). Very recently, ultragraph C*-algebras were connected with the symbolic dynamics of shift spaces over infinite alphabets (see \cite{GRU2}) and ultragraphs were the key object behind a new proposal for the generalization of a shift of finite type to the infinite alphabet case (see \cite{GRU1}).

Due to the exposed above it is natural to study the algebraic analogue of an ultragraph C*-algebra. The formalization of the definition of the algebra was given in \cite{leavittultragraph}, along with a study of the algebra ideals and a proof of a Cuntz-Krieger uniqueness type theorem. Furthermore, it was show in \cite{leavittultragraph} that the class of ultragraph path algebras is strictly larger than the class of Leavitt path algebras. This raises the question of which results about Leavitt path algebras can be generalized to ultragraph path algebras, and whether results from the C*-algebraic setting can be proved in the algebraic level. Our work is a first step in this direction. Building from ideas in \cite{GR1}, where Leavitt path algebras are realized as partial skew group rings, we realize ultragraph path algebras as partial skew group rings. This is also the algebraic version of the characterization of ultragraph C*-algebras as partial crossed products given in \cite{GRU1} (notice that the algebraic version we present is more general than the C*-algebraic version, since the later is valid for ultragraphs with no sinks that satisfy Condition~(RFUM)). 

The theory of partial skew group rings has been in constant development recently, see for example \cite{Gonc, johandanieldanilo} where simplicity criteria are described, and \cite{OinertChain} where chain conditions are studied. In our case we use partial skew ring theory to characterize artinian ultragraph path algebras and give simplicity criteria for these algebras.

Given an ultragraph $\GG$, we realize the associated path algebra as a partial skew group ring in Section~\ref{UPSGR}. For this we consider the free group on the edges of $\GG$. In the graph case (see \cite{GR1}), the free group of edges acts on a subspace of the functions in a set $X$, where $X$ is the set of infinite paths union with finite paths ending in a sink (a vertex that emits no edges). In the ultragraph setting, a finite path is a pair $(\alpha, A)$, where $\alpha= e_1 \ldots e_n$ is a sequence of edges such that $s(e_{i+1})\in r(e_i)$, and $A$ is a subset of $r(e_n)$. To find the correct set $X$ is a key step in our construction. For ultragraphs the set $X$ is formed by the infinite sequences, finite sequences $(\alpha,A)$ such that $A$ contains a sink, and sequences of length zero of the form $(v,v)$ where $v$ is a sink. After defining the set $X$ we proceed with the definition of the partial action and set up the ground to prove Theorem~\ref{isom}, which gives the isomorphism between the partial skew group ring and the ultragraph path algebra. 

In light of Theorem~\ref{isom} we use the results in \cite{johandanieldanilo} to characterize simplicity of ultragraph path algebras in Section 4. As it is the case with Leavitt and graph C*-algebras, the criteria for simplicity we obtain coincides with the one for ultragraph C*-algebras (the later is given in \cite{TomSimple}). More precisely, we show that (when $R$ is a field) the ultragraph Leavitt path algebra is simple if, and only if, $\mathcal{G}$ satisfies Condition~$(L)$ and the unique saturated and hereditary subcollections of $\mathcal{G}^0$ are $\emptyset$ and $\mathcal{G}^0$ (this is Theorem~\ref{simlicitydescribed}). We remark that, using the tools developed in this section, we provide a new proof of the Cuntz-Krieger Uniqueness Theorem for Leavitt path algebras of ultragraphs (Corollary~\ref{ck}). We end the paper in Section~5, where we apply the results of \cite{OinertChain} to characterize artinian ultragraph path algebras.

\section{Ultragraphs and partial skew group rings}

Ultragraph C*-algebras were introduced by Tomforde in \cite{leavittultragraph}. Here we recall the main definitions and relevant results.

\begin{definicao}\label{def of ultragraph}
An \emph{ultragraph} is a quadruple $\mathcal{G}=(G^0, \mathcal{G}^1, r,s)$ consisting of two countable sets $G^0, \mathcal{G}^1$, a map $s:\mathcal{G}^1 \to G^0$, and a map $r:\mathcal{G}^1 \to P(G^0)\setminus \{\emptyset\}$, where $P(G^0)$ stands for the power set of $G^0$.
\end{definicao}

\begin{definicao}\label{def of mathcal{G}^0}
Let $\mathcal{G}$ be an ultragraph. Define $\mathcal{G}^0$ to be the smallest subset of $P(G^0)$ that contains $\{v\}$ for all $v\in G^0$, contains $r(e)$ for all $e\in \mathcal{G}^1$, and is closed under finite unions and non-empty finite intersections.
\end{definicao}


\begin{definicao}\label{def of ultragraph algebra}
Let $\mathcal{G}$ be an ultragraph and $R$ be a unital commutative ring. The Leavitt path algebra of $\mathcal{G}$, denoted by $L_R(\mathcal{G})$ is the universal $R$ with generators $\{s_e,s_e^*:e\in \mathcal{G}^1\}\cup\{p_A:A\in \mathcal{G}^0\}$ and relations
\begin{enumerate}
\item $p_\emptyset=0,  p_Ap_B=p_{A\cap B},  p_{A\cup B}=p_A+p_B-p_{A\cap B}$, for all $A,B\in \mathcal{G}^0$;
\item $p_{s(e)}s_e=s_ep_{r(e)}=s_e$ and $p_{r(e)}s_e^*=s_e^*p_{s(e)}=s_e^*$ for each $e\in \mathcal{G}^1$
\item $s_e^*s_f=\delta_{e,f}p_{r(e)}$ for all $e,f\in \mathcal{G}$
\item $p_v=\sum\limits_{s(e)=v}s_es_e^*$ whenever $0<\vert s^{-1}(v)\vert< \infty$.
\end{enumerate}
\end{definicao}
 
Before we proceed we quickly remind the reader the definition of a partial action: A partial action of a group $G$ on a set $\Omega$ is a pair $\alpha= (\{D_{t}\}_{t\in G}, \ \{\alpha_{t}\}_{t\in G})$, where for each $t\in G$, $D_{t}$ is a subset of $\Omega$ and $\alpha_{t}:D_{t^{-1}} \rightarrow \Delta_{t}$ is a bijection such that $D_{e} = \Omega$, $\alpha_{e}$ is the identity in $\Omega$, $\alpha_{t}(D_{t^{-1}} \cap D_{s})=D_{t} \cap D_{ts}$ and $\alpha_{t}(\alpha_{s}(x))=\alpha_{ts}(x),$ for all $x \in D_{s^{-1}} \cap D_{s^{-1} t^{-1}}.$ In case $\Omega$ is an algebra or a ring then the subsets $D_t$ should also be ideals and the maps $\alpha_t$ should be isomorphisms. 

Associated to a partial action of a group $G$ in a ring $A$ the partial skew group ring, denoted by $A\rtimes_{\alpha} G$, is defined as the set of all finite formal sums $\sum_{t \in G} a_t\delta_t$, where for all $t \in G$, $a_t \in D_t$ and $\delta_t$ is a symbol. Addition is defined component-wise and multiplication is determined by $(a_t\delta_t)(b_s\delta_s) = \alpha_t(\alpha_{-t}(a_t)b_s)\delta_{ts}$

\section{Ultragraph path algebra as a partial skew group ring}\label{UPSGR}

Let $\mathcal{G}$ be an ultragraph. A finite path is either an element of $\mathcal{G}^0$ or a sequence of edges $e_1...e_n$, with length $|e_1...e_n|=n$, and such that $s(e_{i+1})\in r(e_i)$ for each $i\in \{0,...,n-1\}$. An infinite path is a sequence $e_1e_2e_3...$, with length $|e_1e_2...|=\infty$, such that $s(e_{i+1})\in r(e_i)$ for each $i\geq 0$. The set of finite paths in $\mathcal{G}$ is denoted by $\mathcal{G}^*$, and the set of infinite paths in $\mathcal{G}$ is denoted by $\mathfrak{p}^\infty$. We extend the source and range maps as follows: $r(\alpha)=r(\alpha_{|\alpha|})$, $s(\alpha)=s(\alpha_1)$ for $\alpha\in \mathcal{G}^*$ with $0<|\alpha|<\infty$, $s(\alpha)=s(\alpha_1)$ for each $\alpha\in \mathfrak{p}^\infty$, and $r(A)=A=s(A)$ for each $A\in \mathcal{G}^0$. An element $v\in G^0$ is a sink if $s^{-1}(v) = \emptyset$, and we denote the set of sinks in $G^0$ by $G^0_s$. We say that $A\in \mathcal{G}^0$ is a sink if each vertex in $A$ is a sink. 

Define the set $$X=\mathfrak{p}^\infty\cup\{(\alpha,v): \alpha \in \mathcal{G}^*, |\alpha|\geq 1, v \in G^0_s \cap r(\alpha) \}\cup\{(v,v): v\in G^0_s \}  .$$

\begin{remark}
Notice that given a vertex $v$, the element $(v,v)\in X$ if, and only if, $v$ is a sink. 
\end{remark}

\begin{definicao}
For an element $(\alpha,v)\in X$ we define the range and source maps by $r(\alpha,v)=v$ and $s(\alpha,v)=s(\alpha)$. In particular, for a sink $v$, $s(v,v)=v=r(v,v)$. We also extend the length map to the elements $(\alpha,v)$ by defining $|(\alpha,v)|:=|\alpha|$.

\end{definicao}

Next we setup some notation necessary to define the desired partial action.
Let $\F$ be the free group generated by $\mathcal{G}^1$, and denote by $0$ the neutral element of $\F$. Let $W\subseteq \F$ be the set $$W=\{a_1...a_n\in \F: a_i\in \mathcal{G}^1\,\, \forall i\text{ and } s(a_{i+1})\in r(a_i) \forall i\in \{0,...,n-1\}\}.$$ 

\begin{remark} The set $W$ is the same as the set of elements of $\mathcal{G}^*$ with positive length.
\end{remark}

\begin{notation}
 Given an element $a\in W$, with length $|a|$, and an element $x\in X$, we use the notation $x_1...x_{|a|}=a$ to mean that $\alpha_1...\alpha_{|a|}=a_1...a_{|a|}$, if $x=(\alpha,v)\in X$ with $|x|<\infty$, and  $\alpha_1...\alpha_{|a|}=a_1...a_{|a|}$ if $x=\alpha_1 \alpha_2 \ldots$ with $|x|=\infty$.
\end{notation}

Now we define the following sets:
\begin{itemize}
\item for $a\in W$, let $X_a=\{x\in X:x_1..x_{|a|}=a\}$;

\item for $b\in W$, let $X_{b^{-1}}=\{x\in X:s(x)\in r(b)\}$;
\item for $a,b\in W$ with $r(a)\cap r(b)\neq \emptyset$, let $$X_{ab^{-1}}=\left \{x\in X:|x|>|a|, \,\,\,  x_1...x_{|a|}=a \text{ and }s(x_{|a|+1})\in r(b)\cap r(a)\right\}\bigcup$$ $$\bigcup\left\{(a,v) \in X:v\in r(a)\cap r(b)\right\};$$
\item for the neutral element $0$ of $\F$, let $X_0=X$;
\item for all the other elements $c$ of $\F$, let $X_c=\emptyset$.
\end{itemize}

Define, for each $A\in \mathcal{G}^0$ and $b\in W$, the sets 
$$X_A=\{x\in X:s(x)\in A \}$$ and 
$$X_{bA}=\{x\in X_b:|x|>|b| \text{ and } s(x_{|b|+1})\in A\}\cup\{(b,v)\in X_b:v\in A\}.$$

\begin{remark} Notice that for each $a,b\in W$, it holds that $X_{ab^{-1}}=X_{a(r(b))}=X_{a(r(a)\cap r(b))}$ and $X_{r(b)}=X_{b^{-1}}$. Moreover, for an element $b\in W$ and $u\in r(b)$ a sink, it holds that $X_{b\{u\}}=\{(b,u)\}$.
\end{remark}

The following lemma follows from the definitions of the sets $X_c$ and $X_A$, for $c\in \F$ and $A\in \mathcal{G}$, and its proof is left to the reader.

\begin{lema}\label{subsetsintersection} Let $a,b,c,d\in W$ and $A,B\in \mathcal{G}^0$. Then:
\begin{enumerate}
\item $X_a\cap X_b=\left\{\begin{array}{lll}X_a & \text{ if } a=b\xi \text{ for some }\xi\in W\cup\{0\},\\
\emptyset & \text{ if } a_i\neq b_i \text{ for some }i,\\
X_b & \text{ if } b=a\xi \text{ for some } \xi\in W. \end{array} \right . $

\item $X_a\cap X_{c^{-1}}=\left\{\begin{array}{ll}X_a & \text{ if }s(a)\in r(c),\\
\emptyset & \text{otherwise}.\end{array} \right .$

\item $X_a\cap X_{bc^{-1}}=\left\{\begin{array}{ll}X_a & \text{ if } a=b\xi \text{ for some }\xi\in W \text{ and }s(\xi)\in r(c), \\
X_{bc^{-1}} & \text{ if } b=a\xi \text{ for some }\xi\in W\cup \{0\},\\
\emptyset & \text{otherwise}.

\end{array} \right .$

\item\label{item4subset} $X_{ab^{-1}}\cap X_{cd^{-1}}=\left\{\begin{array}{ll}X_{ab^{-1}} & \text{ if } a=c\xi \text{ for some }\xi \in W \text{ and }s(\xi)\in r(d), \\

X_{cd^{-1}} & \text{ if } c=a\xi \text{ for some }\xi \in W \text{ and }s(\xi)\in r(b), \\
X_{a(r(b)\cap r(d)) } & \text{ if } a=c,\\

\emptyset & \text{otherwise}.
\end{array} \right .$

\item\label{item5} $X_A\cap X_a=\left\{\begin{array}{ll}X_a & \text{ if } s(a)\in A,\\
\emptyset & \text{otherwise}.
\end{array} \right .$

\item\label{XAXab} $X_A\cap X_{ab^{-1}}=\left\{\begin{array}{ll}X_{ab^{-1}} & \text{ if } s(a)\in A,\\
\emptyset & \text{otherwise}.
\end{array} \right .$
\item\label{item7} $X_A\cap X_B=X_{A\cap B}$ and $X_A\cup X_B=X_{A\cup B}$.

\item\label{XbAXc} $X_{bA}\cap X_c=\left\{\begin{array}{ll} X_{bA} & \text{ if } b=c\xi \text{ for some } \xi\in W\cup \{0\},\\
X_c & \text{ if } c=b\xi \text{ for some } \xi\in W \text{ 
and } s(\xi)\in A,\\ 
\emptyset & \text{ otherwise. }
\end{array}\right .$

\item\label{item9} $X_{bA}\cap X_{cd^{-1}}=\left\{\begin{array}{ll}
X_{bA} & \text{ if } b=c\xi \text{ for some } \xi\in W \text{ and } s(\xi)\in r(d),\\
X_{cd^{-1}} & \text{ if } c=b\xi \text{ for some } \xi\in W \text{ and } s(\xi)\in A,\\
X_{b(A\cap r(d))} & \text{ if }b=c,\\
\emptyset & \text{ otherwise. }
\end{array}\right .$

\item $X_{bA}\cap X_{cB}=\left\{\begin{array}{ll}
X_{bA} & \text{ if } b=c\xi \text{ for some } \xi\in W \text{ and } s(\xi)\in B,\\
X_{cB} & \text{ if } c=b\xi \text{ for some } \xi\in W \text{ and } s(\xi)\in A,\\
X_{b(A\cap B)} & \text{ if } b=c,\\
\emptyset & \text{ otherwise. }

\end{array}\right .$

\end{enumerate}
\end{lema}

Our aim is to get a partial action from $\F$ on $X$. With this in mind, define the following bijective maps:
\begin{itemize}
\item for $a\in W$ define $\theta_a:X_{a^{-1}}\rightarrow X_a$ by $$\theta_a(x)=\left\{\begin{array}{ll}
ax & \text{ if } |x|=\infty,\\
(a\alpha,v) & \text{ if } x=(\alpha,v),\\
(a,v) & \text{ if } x=(v,v);

\end{array} \right .$$

\item for $a\in W$ define $\theta_a^{-1}:X_a\rightarrow X_{a^{-1}}$ as being the inverse of $\theta_a$;
\item for $a,b\in W$ define $\theta_{ab^{-1}}:X_{ba^{-1}}\rightarrow X_{ab^{-1}}$ by 

$$\theta_{ab^{-1}}(x)=\left\{\begin{array}{ll}
ay & \text{ if } |x|=\infty \text{ and } x=by,\\
(a\alpha,v) & \text{ if } x=(b\alpha,v),\\
(a,v) & \text{ if } x=(b,v);
\end{array} \right .$$

\item for the neutral element $0\in \F$ define $\theta_0:X_0\rightarrow X_0$ as the identity map;

\item for all the other elements $c$ of $\F$ define $\theta_c:X_{c^{-1}}\rightarrow X_c$ as the empty map.

\end{itemize}

\begin{remark}\label{XaB} Notice that $$X_{bA}=\{x\in X_b;\theta_{b^{-1}}(x)\in X_A)\}=\{x\in X_b;\theta_{b^{-1}}(x)\in X_A\cap X_{b^{-1}}\}=$$ $$=\theta_b(X_A\cap X_{b^{-1}}), $$ that is,  $X_{bA}=\theta_b(X_A\cap X_{b^{-1}}).$
\end{remark}

It is straightforward to check that $(\{\theta_t\}_{t\in \F}, \{X_t\}_{t\in \F})$ is a partial action of $\F$ on $X$, that is, $X_e=X$, $\theta_e=Id_x$, $\theta_c(X_{c^{-1}}\cap X_t)=X_{ct}\cap X_c$ and $\theta_c\circ\theta_t=\theta_{ct}$ in $X_{t^{-1}}\cap X_{t^{-1}c^{-1}}$. Define for each $c\in \F$ the set $F(X_c)$ of all the functions from $X_c$ to the commutative unital ring $R$. Notice that each $F(X_c)$ is an $R$-algebra, with pointwise sum and product. For the neutral element $0\in \F$ we denote the set $F(X_0)$ simply by $F(X)$. Each $F(X_c)$ is an ideal of the $R$-algebra $F(X)$. Now, for each $c\in \F$ define the $R$-isomorphism $$\beta_c:F(X_c^{-1})\rightarrow F(X_c)$$ by $\beta_c(f)=f\circ \theta_{c^{-1}}$, whose inverse is the isomorphism $\beta_{c^{-1}}$. So, we get a partial action $(\{\beta_c\}_{c\in \F}, \{F(X_c)\}_{c\in \F})$ from $\F$ to the $R$-algebra $F(X)$.

To get the desired partial action we need to restrict the partial action $\beta$ to the $R$-subalgebra $D$ of $F(X)$ generated by all the finite sums of all the finite products of the characteristic maps $\{1_{X_A}\}_{A\in \mathcal{G}^0}$, $\{1_{bA}\}_{b\in W, A\in \mathcal{G}^0}$ and $\{1_{X_c}\}_{c\in \F}$. We also define, for each $t\in \F$ the ideals $D_t$ of $D$, as being all the finite sums of finite products of the characteristic maps $\{1_{X_t}1_{X_A}\}_{A\in \mathcal{G}^0}$, $\{1_{X_t}1_{bA}\}_{b\in W, A\in \mathcal{G}^0}$ and $\{1_{X_t}1_{X_c}\}_{c\in \F}$.

\begin{remark}\label{spanDt} From now on we will use the notation $1_A$, $1_{bA}$ and $1_t$ instead of $1_{X_A}$, $1_{X_{bA}}$ and $1_{X_t}$, for $A\in \mathcal{G}^0$, $b\in W$ and $t\in \F$. It follows directly from Lemma \ref{subsetsintersection} that 
$$D=\text{span}\{1_A, 1_c, 1_{bA}:A\in \mathcal{G}^0, c\in \F\setminus\{0\}, b\in W\}, $$ and that for each $t\in \F$, 
$$D_t=\text{span}\{1_t1_A, 1_t1_c, 1_t1_{bA}:A\in \mathcal{G}^0, c\in \F, b\in W\}, $$ where ``span'' means linear span. 
\end{remark}

Our aim is to restrict the partial action $\beta$ to the ideals $\{D_t\}_{t\in \F}$ of $D$. The next proposition tells us that $\beta_t(D_{t^{-1}})=D_t$ for each $t\in \F$.

\begin{proposicao}\label{actionbeta} 
\begin{enumerate}
\item\label{item1actionbeta}  For $t,c\in \F$ it holds that $\beta_c(1_{c^{-1}}1_t)=1_c1_{ct}$.
\item\label{item2actionbeta} For $b\in W$ and $A\in \mathcal{G}^0$ we get $\beta_b(1_b^{-1}1_A)=1_b1_{bA}$.

\item For $t=ab^{-1}$ with $b\in W$ and $a\in W\cup\{0\}$, and $A\in \mathcal{G}^0$, we get $$\beta_t(1_{t^{-1}}1_A)=\left\{\begin{array}{ll}1_t & \text{ if } s(b)\in A,\\
0 & \text{ otherwise. }

\end{array}\right .$$

\item For $b,c\in W$ and $A\in \mathcal{G}^0$ it holds that
$$\beta_c(1_{c^{-1}}1_{bA})=\left\{\begin{array}{ll}
1_c1_{cbA} & \text{ if } s(b)\in r(c),\\
0 & \text{ otherwise. }
\end{array}\right .$$

\item For $b,c,d\in W$, and $A\in \mathcal{G}^0$, we get
$$\beta_{dc^{-1}}(1_{cd^{-1}}1_{bA})=\left\{\begin{array}{ll}
1_{dc^{-1}} 1_{d\xi A} = 1_{d\xi A} & \text{ if } b=c\xi \text{ for some } \xi\in W \text{ and } s(\xi) \in r(d),\\
1_{dc^{-1}} & \text{ if } c=b\xi \text{ for some } \xi \in \text{ and } s(\xi)\in A ,\\
1_{dc^{-1}}1_{dA} = 1_{d(r(c)\cap A)} & \text{ if } b=c,\\
0 & \text{ otherwise. }
\end{array}\right .$$

\item For $a,b,c\in W$ and $A\in \mathcal{G}^0$, we get
$$\beta_{c^{-1}}(1_c1_{bA})=\left\{\begin{array}{ll}
1_{c^{-1}}1_{\xi A} = 1_{\xi A} & \text{ if } b=c\xi \text{ for some } \xi \in W\cup\{0\}\\ 
1_{c^{-1}} & \text{ if } c=b\xi \text{ for some } \xi \in W \text{ and } s(\xi)\in A\\
0 & \text{ otherwise }

\end{array}\right .$$

\end{enumerate}
\end{proposicao}

\demo The first item follows from the fact that $\theta_c(X_{c^{-1}}\cap X_t)=X_c\cap X_{ct}$, since $\beta_c(1_{c^{-1}}1_t)(x)=[\theta_{c^{-1}}(x)\in X_{c^{-1}}\cap X_t]=[x\in \theta_c(X_{c^{-1}}\cap X_t)]=[x\in (X_c\cap X_{ct})]=1_c(x)1_{tc}(x)$.

To see that the second item holds, note that $\beta_b(1_{b^{-1}}1_A)(x)=[\theta_{b^{-1}}(x)\in (X_{b^{-1}}\cap X_A)]=[x\in \theta_b(X_{b^{-1}}\cap X_A)]=[x\in X_{bA}]=1_{bA}(x)$, where the second to last equality follows from Remark \ref{XaB}.
The third item follows from Item \ref{XAXab} of Lemma~\ref{subsetsintersection}.

To see that Item 4. holds note that, for $x\in X_c$, 
$$\beta_c(1_{bA}1_{c^{-1}})(x)=[\theta_{c^{-1}}(x)\in X_{bA}\cap X_{c^{-1}}]=[x\in \theta_c(X_{bA}\cap X_{c^{-1}})]=$$ $$=[x\in X_{cbA}\cap X_c]=1_c(x)1_{cbA}(x),$$ and for $x\notin X_c$,  $\beta_c(1_{bA}1_{c^{-1}})(x)=0=1_c(x)1_{cbA}(x).$

Item 5 follows from Item \ref{item9} of 
Lemma~\ref{subsetsintersection}, and the last item follows from Item \ref{XbAXc} of the same Lemma.
\fim

By the previous proposition we get that, for each $t\in \F$, $\beta_t(D_{t^{-1}})\subseteq D_t$ and, consequently, $\beta_t(D_{t^{-1}})= D_t$ for each $t\in \F$. So we may consider the restriction of the partial action $\beta$ to the subsetes $\{D_t\}_{t\in \F}$ of $D$. We denote this restriction also by $\beta$, and so we get a partial action $(\{\beta_t\}_{t\in \F}, \{D_t\}_{t\in \F})$ of $\F$ in $D$.
 Now we are ready to prove the following theorem.

\begin{teorema}\label{isom} Let $\mathcal{G}$ be an ultragraph, $R$ be an unital commutative ring, and let $L_R(\mathcal{G})$ be the Leavitt path algebra of $\mathcal{G}$. Then there exists an $R$-isomorphism $\phi:L_R(\mathcal{G})\rightarrow D\rtimes_\beta \F$ such that $\phi(p_A)=1_A\delta_0$, $\phi(s_e^*)=1_{{e^{-1}}}\delta_{e^{-1}}$ and $\phi(s_e)=1_e\delta_e$ for each $A\in \mathcal{G}^0$ and $e\in \mathcal{G}^1$.
\end{teorema}

\demo First we show that the sets $\{1_A\delta_0\}_{A\in \mathcal{G}^0}$ and $\{1_e\delta_e, 1_{e^{-1}}\delta_{e^{-1}}\}_{e\in \mathcal{G}^1}$ satisfies the relations which define the algebra $L_R(\mathcal{G})$. 

The first relation of Definition 
\ref{def of ultragraph algebra} follows from Item~\ref{item7} of Lemma~\ref{subsetsintersection}. To verify the second relation, let $e\in \mathcal{G}^1$, and note that
$1_{s(e)}\delta_0 1_e\delta_e=1_{s(e)}1_e\delta_e=1_e\delta_e$ and $1_e\delta_e1_{r(e)}\delta_0=\beta_e(\beta_{e^{-1}}(1_e)1_{r(e)})\delta_e=\beta_e(1_{e^{-1}}1_{r(e)})\delta_e=1_e1_{er(e)}\delta_e=1_e\delta_e$, where the second to last equality follows from Item \ref{item2actionbeta} of Proposition \ref{actionbeta}.
Moreover, $1_{r(e)}\delta_01_{e^{-1}}\delta_{e^{-1}}=1_{e^{-1}}\delta_{e^{-1}}$ and $1_{e^{-1}}\delta_{e^{-1}}1_{s(e)}\delta_e=\beta_{e^{-1}}(\beta_e(1_{e^{-1}})1_{s(e)})\delta_{e^{-1}}=\beta_{e^{-1}}(1_e1_{s(e)})\delta_{e^{-1}}=\beta_{e^{-1}}(1_e)\delta_{e^{-1}}=1_{e^{-1}}\delta_{e^{-1}}$. Next we verify the third relation. Let $e,f\in \mathcal{G}^1$. Then $$1_{e^{-1}}\delta_{e^{-1}}1_f\delta_f=\beta_{e^{-1}}(1_e1_f)\delta_{e^{-1}f}.$$ If $e\neq f$ then $1_e1_f=0$ and if $e=f$ then $\beta_{e^{-1}}(1_e1_f)\delta_{e^{-1}f}=\beta_{e^{-1}}(1_e)\delta_0=1_{e^{-1}}\delta_0=1_{r(e)}\delta_0$. To verify the last relation of Definition~\ref{def of ultragraph algebra}, note first that $1_e\delta_e1_{e^{-1}}\delta_{e^{-1}}=1_e\delta_0$, for each edge $e$. Now, let $v$ be an vertex such that $0<|s^{-1}(v)|<\infty$. Then $X_v=\bigcup\limits_{e\in s^{-1}(v)} X_e$, from where $1_v=\sum\limits_{e\in s^{-1}(v)}1_e$, and so $$\sum\limits_{e\in s^{-1}(v)}1_e\delta_e1_{e^{-1}}\delta_{e^{-1}}=\sum\limits_{e\in s^{-1}(v)}1_e\delta_0=1_v\delta_0.$$ 

So, by the universality of $L_R(\mathcal{G})$, there exists an $R$-homomorphism $\phi:L_R(\mathcal{G})\rightarrow D\rtimes_\beta\F$ such that $\phi(p_A)=1_A\delta_0$, $\phi(s_e)=1_e\delta_e$ and $\phi(s_e^*)=1_{e^{-1}}\delta_{e^{-1}}$ for each $A\in \mathcal{G}^0$ and each edge $e$.

Now we prove that $\phi$ is surjective. 

For each $a=a_1...a_{|a|}\in W$ and $d=d_1...d_{|d|}\in W$ we use the notations $\phi(s_a)$, $\phi(s_d^*)$ and $\phi(s_as_d^*)$ to denote the elements $\phi(s_{a_1})...\phi(s_{a_{|a|}})$, $\phi(s_{d_{|d|}}^*)...\phi(s_{d_1}^*)$ and $\phi(s_{a_1})...\phi(s_{a_{|a|}})\phi(s_{d_{|d|}}^*)...\phi(s_{d_1}^*)$ respectively.

{\it Claim 1: For each $a,d\in W$ it holds that $\phi(s_a)\phi(s_a^*)=1_a\delta_0$, $\phi(s_d^*)\phi(s_d)=1_{d^{-1}}\delta_0$ and $\phi(s_as_d^*)\phi(s_ds_a^*)=1_{ad^{-1}}\delta_0$.}  

The equalities $\phi(s_a)\phi(s_a^*)=1_a\delta_0$ and $\phi(s_d^*)\phi(s_d)=1_{d^{-1}}\delta_0$ follow by induction on the length of $a$ and $d$ and from the first item of Proposition \ref{actionbeta}. To prove the other equality write $a=eg$, where $|e|=1$ and $|g|=|a|-1$, and suppose by inductive arguments that $\phi(s_g)\phi(s_d^*)\phi(s_d)\phi(s_g)^*=1_{gd^{-1}}\delta_0$. Then $$\phi(s_as_d^*)\phi(s_ds_a^*)=\phi(s_e)\phi(s_g)\phi(s_d^*)\phi(s_d)\phi(s_g ^*)\phi(s_e)=\phi(s_e)1_{gd^{-1}}\phi(s_e)=$$ $$=1_e\delta_e1_{gd^{-1}}1_{e^{-1}}\delta_{e^{-1}}=\alpha_e(1_{e^{-1}}1_{gd^{-1}})\delta_0=1_e1_{egd^{-1}}\delta_0=1_{ad^{-1}}\delta_0,$$ where the second to last equality follows from the first item of Proposition \ref{actionbeta}. So, Claim 1 is proved.

{\it Claim 2: For each $b\in W$, and $A\in \mathcal{G}^0$, it holds that $\phi(s_b)\phi(p_A)\phi(s_b^*)=1_{bA}\delta_0$.}

For $|b|=1$ note that $\phi(s_b)\phi(p_A)\phi(s_b^*)=\beta_b(1_{b^{-1}}1_A)\delta_0=1_b1_{bA}\delta_0=1_{b_A}\delta_0$, where the second to last equality follows from Item \ref{item2actionbeta} of Proposition \ref{actionbeta}. Now, for $|b|>1$, write $b=ed$ with $|e|=1$ and $|d|=|b|-1$. By inductive arguments we get that
$$\phi(s_b)\phi(p_A)\phi(s_b^*)=\phi(s_e)\phi(s_d)\phi(p_A)\phi(s_d^*)\phi(s_e^*)=$$ $$=\phi(s_e) 1_{dA}\delta_0\phi(s_e^*)=\beta_e(1_{e^{-1}}1_{dA})\delta_0=1_e1_{edA}\delta_0=1_{bA}\delta_0,$$ where the second to last equality follows by similar arguments to the ones used in the proof of Item \ref{item2actionbeta} of Proposition \ref{actionbeta}. So, Claim 2 is proved.

By Remark~\ref{spanDt}, to prove that $\phi$ is surjective,  it is enough to prove that $$\{ 1_A\delta_0, 1_c\delta_0, 1_{bA}\delta_0: A\in \mathcal{G}^0, c\in \F\setminus \{0\}, b\in W\}\subseteq Im(\phi)$$ and, for each $t\in \F$,

$$\{1_t1_A\delta_t, 1_t1_c\delta_t, 1_t1_{bA}\delta_t:A\in\mathcal{G}^0, c\in \F, b\in W\}\subseteq Im(\phi).$$

{\it Claim 3: $\{ 1_A\delta_0, 1_c\delta_0, 1_{bA}\delta_0: A\in \mathcal{G}^0, c\in \F\setminus\{0\}, b\in W\}\subseteq Im(\phi)$.}

Recall that for each $A\in \mathcal{G}^0$, $\phi(p_A)=1_A\delta_0$. Moreover, for $c\in \F\setminus \{0\}$ with $c=ad^{-1}$, where $a,d\in W\cup\{0\}$, we get by Claim 1 that $1_c\delta_0\in Im(\phi)$ (for all the other $c\in \F\setminus \{0\}$ we also have $1_c\delta_0\in Im(\phi)$, since $1_c=0$). To finish notice that, by Claim 2, we get that $1_{bA}\delta_0\in Im(\phi)$, for each $b\in W$ and $A\in \mathcal{G}^0$. So, Claim 3 is proved.

{\it Claim 4: For each $t\in \F\setminus\{0\}$,  $$\{ 1_t1_A\delta_t, 1_t1_c\delta_t, 1_t1_{bA}\delta_t: A\in \mathcal{G}^0, c\in \F\setminus\{0\}, b\in W\}\subseteq Im(\phi).$$}
First, for $e\in W$, with $|e|=1$, recall that $1_e\delta_e=\phi(e)$. Now, let $c\in W$ with $|c|>1$, write $c=ed$ with $|e|=1$, and suppose (by inductive arguments on $|c|$) that $\phi(d)=1_d\delta_d$. Then $$\phi(s_c)=\phi(s_e)\phi(s_d)=1_e\delta_e 1_d \delta_d=\beta_e(1_{e^{-1}}1_d)\delta_{ed}=1_{e1_{ed}}\delta_{ed}=1_c\delta_c,$$ where the second to last equality follows from Item \ref{item1actionbeta} of Proposition \ref{actionbeta}. Analogously we get that $\phi(s_d^*)=1_{d^{-1}}\delta_{d^{-1}}$ for each $d\in W$. Now, for $c,d\in W$, $$\phi(s_c)\phi(s_d^*)=1_c\delta_c1_{d^{-1}}\delta_{d^{1}}=\beta_c(1_{c^{-1}}1_{d^{-1}})\delta_{cd^{-1}}=1_c 1_{cd^{-1}}\delta_{cd^{-1}}=1_{cd^{-1}}\delta_{cd^{-1}},$$ where, again, the second to last equality follows from Item \ref{item1actionbeta} of Proposition \ref{actionbeta}. So we get $1_t\delta_t\in Im(\phi)$ for each $t\in \F\setminus\{0\}$.

Now, for $t,c\in \F\setminus\{0\}$, $b\in W$ and $A\in \mathcal{G}^0$, note that $1_t1_{bA}\delta_t=1_{bA}\delta_01_t\delta_t\in Im(\phi)$, and similarly one shows that $1_t1_A\delta_t, 1_t1_c\delta_t\in Im(\phi)$. So, we get that $\phi$ is surjective.

It remains to show that $\phi$ is injective. To prove this we will use the graded uniqueness theorem, see \cite[Theorem 3.2]{leavittultragraph}. For each integer number $n$ define 
$$F_n=span\{f_{ab^{-1}}\delta_{ab^{-1}}: f_{ab^{-1}}\in D_{ab^{-1}}, \,\, a,b\in W\cup\{0\} \text{ and } |a|-|b|=n\}.$$ Note that $D\rtimes_\beta \F$ is $\Z$-graded by the gradation $\{F_n\}_{n\in \Z}$. Moreover, $L_R(\mathcal{G})$ is a $\Z$-graded ring with the grading $L_R(\mathcal{G})_n=span\{s_ap_As_b^*:a,b\in \mathcal{G}^*, A\in \mathcal{G}^0\}$ introduced in 
\cite{leavittultragraph}. It is easy to see that $\phi$ is a graded ring homomorphism. Since $X_A\neq \emptyset$ then $\phi(\tau p_A)=\tau 1_{A}\neq 0$, for each $A\in \mathcal{G}^0$ and $\tau\in R\setminus\{0\}$. It follows from \cite[Theorem 3.2]{leavittultragraph} that $\phi$ is injective and hence an isomorphism.

\fim

\section{Simplicity and maximal commutativity}

In  this section we use the realization of ultragraph Leavitt path algebras as partial skew group rings to describe simplicity criteria for these algebras. Recall that from \cite[Theorem~2.3]{johandanieldanilo}, the algebra $D\rtimes_\beta\F$ is simple if, and only if, $D$ is $\F$-simple and $D\delta_0$ is maximal commutative in $D\rtimes_\beta\F$. Aiming at the simplicity criteria given for ultragraph C*-algebras in \cite{TomSimple} we will characterize maximal commutativity in terms of Condition~$(L)$ and $\F$ simplicity in terms of hereditary and saturated subcollections $\GG^0$.

Recall that a cycle in an ultragraph $\mathcal{G}$ is a path $\alpha=e_1...e_{|\alpha|}$, with $|\alpha|\geq 1$ and $s(\alpha)\in r(\alpha)$, and an exit for $\alpha$ is an edge $e$ with $s(e)=s(e_i)$, for some $i\in \{1,..., |\alpha|\}$ and $e\neq e_i$. The ultragraph $\mathcal{G}$ satisfies Condition~$(L)$ if each cycle $\alpha=e_1...e_{|\alpha|}$ has an exit, or if $r(e_i)$ contains a sink for some $i$.

Before we state our next result we recall the notion of maximal commutativity: The centralizer of a nonempty subset $S$ of a ring $R$, which we denote by $C_R(S)$, is the set of all elements of $R$ that commute with each element of $S$. If $C_R(S)=S$ holds, then $S$ is said to be a \emph{maximal commutative subring} of $R$.

\begin{teorema}\label{maximalcommutative} Let $\mathcal{G}$ be an ultragraph. Then $D\delta_0$ is maximal commutative in $D\rtimes_\beta\F$ if, and only if, $\mathcal{G}$ satisfies condition $(L)$.
\end{teorema}

\demo First suppose that $\mathcal{G}$ satisfies condition $(L)$. Suppose, by contradiction, that there exists $x = \sum a_t \delta_t$, with some $t\neq 0$, that commutes with $a\delta_0$ for all $a \in D$. Then there exists $t\in \F\setminus\{0\}$, and $a_t\in D_t$ with $a_t\neq 0$, such that $a_t\delta_ta_0\delta_0=a_0\delta_0a_t\delta_t$ for each $a_0\in D$. From the last equality we get 
\begin{equation}\label{eq1}
\beta_t(\beta_{t^{-1}}(a_t)a_0)=a_ta_0
\end{equation} 
for each $a_0\in D_0$. Since $a_t\neq 0$ then either $t=a$, $t=b^{-1}$, or $t=ab^{-1}$, with $a,b\in W$. 

Notice that, since $$D_t=span\{1_t1_c, 1_t1_{bA}, 1_t1_A:c\in \F, b\in W, A\in \mathcal{G}^0\}$$ then, for each $\xi\in X_t$ with $|\xi|=\infty$, there exists an $m\in \N$ such that, if $\eta\in X_t$ and $\eta_1\eta_2...\eta_m=\xi_1\xi_2...\xi_m$ then $a_t(\eta)=a_t(\xi)$.

We now divide the proof in three cases.

{\it Case 1: Suppose $t\in W$.}

If we take $a_0=1_{t^{-1}}$ in Equation~(\ref{eq1}) we get that $a_t=a_t 1_{t^{-1}}$. Hence the support of $a_t$ is contained in $X_t\cap X_{t^{-1}}$, and therefore $t$ is a closed path. If we take $a_0=1_t1_{t^{-1}}$ then, from Equation (\ref{eq1}), we have that $\beta_t(\beta_{t^{-1}}(a_t)1_t)=a_t1_t1_{t^{-1}}=a_t$, and from Remark \ref{spanDt} and Proposition \ref{actionbeta}, we get $\beta_t(\beta_{t^{-1}}(a_t)1_t)=a_t1_{tt}$. Therefore $a_t1_{tt}=a_t$. With the same arguments, if we take $a_0=1_{t^2}$ we get $a_t1_{t^3}=a_t$, and inductively we get $a_t1_{t^n}=a_t$ for each $n\in \N$.

Let $\xi \in X_t$ be such that $a_t(\xi)\neq 0$. Then $a_t(\xi)1_{t^n}(\xi)\neq 0$, for each $n\in \N$, and so $|\xi|=\infty$. Let $m\in \N$ be such that if $\eta\in X_t$, and $\eta_1...\eta_m=\xi_1...\xi_m$, then $a_t(\eta)=a_t(\xi)$.

Since $\mathcal{G}$ satisfies condition $(L)$ the closed path $t=t_1...t_{|t|}$ either has an exit or some $r(t_i)$ contains a sink. 

Suppose first that $t$ has an exit, that is, there exists an edge $e$ such that $s(e)\in r(t_i)$, for some $i$ and $e\neq t_{i+1}$. Let $k\in \N$ be such that $k|t|\geq m$ and let $\eta$ be such that $\eta=t^{k}t_1t_2...t_iey$ (for some $y$). Then we get that $0\neq a_t(\xi)=a_t(\eta)=(a_t 1_{t^{k+1}})(\eta)=0$, a contradiction.

Now suppose that $r(t_i)$ contains a sink $v$ for some $i$. Then, again, let $k\in \N$ be such that $k|t|\geq m$, and let $\eta=(t^kt_1t_2...t_i,v)$, which is an element of $X_t$. Then we have that $0\neq a_t(\xi)=(a_t1_{t^{k+1}})(\xi)=(a_t1_{t^{k+1}})(\eta)=0$, which is also a contradiction.

So we conclude that $t\notin W$.

{\it Case 2: $t=d^{-1}$, with $d\in W$.}

From Equation (\ref{eq1}) we get that $\beta_{d^{-1}}(\beta_d(a_{d^{-1}})a_0)=a_{d^{-1}}a_0$ and so $\beta_d(a_{d^{-1}})a_0=\beta_d(a_{d^{-1}}a_0)$. Let $c_d=\beta_d(a_{d^{-1}})$. Then $\beta_{d^{-1}}(c_d)=a_{d^{-1}}$ and so we get the equality $$\beta_d(\beta_{d^{-1}}(c_d)a_0)=c_da_0,$$ for each $a_0\in D_0$. Now, by {\it Case 1}, we get a contradiction and hence it is not possible that $t=d^{-1}$ with $d\in W$.

{\it Case 3: $t=cd^{-1}$ with $c,d\in W$}.

As in {\it Case 1} we get that $a_t=a_t1_{t^n}$ for each $n\in \N$. Hence, since $a_t\neq 0$, we have that $X_{t^n}\neq \emptyset$ for each $n$. Therefore 
either $c=db$ or $d=cb$ with $b\in W$. 

If $c=d b$ then $t^n=db^nd^{-1}$ and so $b$ is a closed path. Let $\xi\in X_t$ with $|\xi|=\infty$ and $a_t(\xi)\neq 0$. Proceeding from this point as {\it Case 1} we get a contradiction. 

If $d=cb$ for some $b\in W$ then we also get a similar contradiction, by considering the equality $\beta_{t^{-1}}(\beta_t(u_{t^{-1}})a_0)=u_{t^{-1}}a_0$ obtained from Equation (\ref{eq1}), where $u_{t^{-1}}=\beta_{t^{-1}}(a_t)$.

So, we proved that if $\mathcal{G}$ satisfies condition $(L)$ then $D$ is maximal commutative in $D\rtimes_\beta \F$. Next we prove the converse.

Suppose that $\mathcal{G}$ does not satisfy condition $(L)$. Then there exist a closed path $t=t_1...t_{|t|}$ in $\mathcal{G}$ such that $t$ has no exit and $r(t_i)$ contains exactly one vertex, for each $t_i$. We show that $1_t\delta_t$ commutes with $D_0\delta_0$. By Remark \ref{spanDt} it is enough to show that $1_t\delta_t$ commutes with $1_c\delta_0$ for each $c\in \F\setminus\{0\}$, and with $1_A\delta_0$ and $1_{bA}\delta_0$ for each $A\in\mathcal{G}^0$ and $b\in W$.

Let $A\in \mathcal{G}^0$. If $r(t)=s(t)\in A$ then $1_A\delta_0 1_t\delta_t=1_t\delta_t=\beta_t(1_{t^{-1}})\delta_t=\beta_t(1_{t^{-1}}1_A)\delta_t=\beta_t(\beta_{t^{-1}}(1_t)1_A)\delta_t=1_t\delta_t 1_A\delta_0$,
and if $s(t)=r(t)\notin A$ then $1_A\delta_0 1_t\delta_t=0=1_t\delta_t1_A\delta_0$.

Now let $A\in \mathcal{G}^0$ and $b\in W$. Note that $1_t\delta_t1_{bA}\delta_0=\beta_t(1_{t^{-1}}1_{bA})\delta_t$ and $1_{bA}\delta_0 1_t\delta_t=1_{bA}1_t\delta_t$. If $s(b)\notin r(t)$ then $1_{t^{-1}}1_{bA}=0=1_t1_{bA}$ and we are done. Suppose that $s(b)\in r(t)$. Then, by Proposition \ref{actionbeta},  $\beta_t(1_{t^{-1}}1_{bA})=1_t1_{tbA}$. So, it remains to show that $1_t1_{tbA}=1_t1_{bA}$. Notice that $X_t=\{\xi\}$, where $\xi$ is the infinite path $\xi=tt...$. Then to verify the desired equality it is enough to show that $\xi\in X_{tbA}$ if, and only if, $\xi\in X_{bA}$.
 Suppose that $\xi\in X_{tbA}$. Then $\xi=tby$, where $y$ is a path such that $s(y)\in A$. Therefore there exists an $n\in \N$ such that $b=t^nt_1...t_i$ for some $i$ and note that $s(y)=r(t_i)$. Hence, $$\xi=tby=tt^{n}t_1...t_iy=t^nt_1...t_it_{i+1}...t_{|t|}t_1...t_iy=bt_{i+1}...t_{|t|}t_1...t_iy.$$
Now note that $bt_{i+1}...t_{|t|}t_1...t_iy\in X_{bA}$, since $s(t_{i+1})=r(t_i)=s(y)\in A$.
Similarly one shows that if $\xi\in X_{bA}$ then $\xi\in X_{tbA}$. So, $1_t1_{tbA}=1_t1_{bA}$.

Finally, we show that $1_t\delta_t 1_c\delta_0=1_c\delta_01_t\delta_t$, for each $c\in \F\setminus \{0\}$. To prove this it is sufficient to show that $\beta_t(1_{t^{-1}}1_c)=1_t1_c$, for each $c\in \F\setminus \{0\}$. By Proposition~\ref{actionbeta} we have that $\beta_t(1_{t^{-1}}1_c)=1_t1_{tc}$, and hence we have to show that $1_t1_{tc}=1_t1_c$. Notice that to prove this last equality it is enough to show that $\xi=tt...$ is an element of $X_{tc}$ if, and only if, $\xi\in X_c$. This follows by arguments similar to the previous case, splitting the proof in cases depending whether $c=a$, $c=b$ or $c=ab^{-1}$ with $a,b\in W$. 
\fim

The next proposition will be useful in the characterization of $\F$ simplicity of $D$.

\begin{proposicao}\label{1vtheorem} Let $x_0\delta_0$ be a non-zero element of $D\delta_0$ and let $I$ be the ideal generated by $x_0\delta_0$ in $ D\rtimes_\beta\F$. Then there there exists a vertex $v\in G^0$, and a non-zero element $h\in R$, such that $(h 1_v)\delta_0\in I$.
\end{proposicao}

\demo First note that by Remark~\ref{spanDt} $$x_0=\sum\limits_{i=1}^m\alpha_i1_{a_ib_i^{-1}}+\sum\limits_{j=1}^n\beta_j1_{e_jA_j}+\sum\limits_{k=1}^p\gamma_k1_{B_k},$$ with $a_i, b_i, e_j\in W$ and $a_ib_i^{-1}\neq 0$, $A_j, B_k\in \mathcal{G}^0$, and $\alpha_i,\beta_j, \gamma_k \in R$. Let $A=\{s(a_i): 1\leq i\leq m\}\cup\{s(e_j):1\leq j\leq n\}\bigcup\limits_{k=1}^p B_k$, which is an element of $\mathcal{G}^0$, and note that $1_Ax_0=x_0$. Let $\xi\in X$ be such that $x_0(\xi)\neq 0$, and let $v=s(\xi)$. Then $v\in A$ and so $1_v(\xi)x_0(\xi)=1_A(\xi)x_0(\xi)=x_0(\xi)\neq 0$. Therefore $1_vx_0\neq 0$. 

If $v$ is a sink then $1_vx_0=\sum\limits_{k=1}^p\gamma_k1_v1_{B_k}=\sum\limits_{k\in\{1...p\}: v\in B_k}\gamma_k1_v=h1_v$. So, $(h1_v)\delta_0 \in I$. 

Now suppose that $v$ is not a sink. Let $M=max\{|a_i|, |e_j|: 1\leq i\leq m, 1\leq j \leq n\}$. Note that since $v$ is not a sink then 

$$X_v=\bigcup\limits_{c\in J}^.X_{c\{u\}}\bigcup\limits_{L}^.X_d,$$ where $c$ and $d$ are all the elements of $W$ such that $s(c)=v$, $|c|<M$ and $u\in r(c)$ is a sink, and $s(d)=v$ and $|d|=M+1$. 

Since $1_vx_0\neq 0$ then $1_{c\{u\}}x_0\neq 0$, for some $c\in J$ and some sink $u\in r(c)$, or $1_dx_0\neq 0$ for some $d\in L$.

Suppose that $1_{c\{u\}}x_0\neq 0$. Note that for each $i\in \{1,...,m\}$, $j\in \{1,...,n\}$ and $k\in \{1,...,p\}$, we have that $1_{c\{u\}}1_{a_ib_i^{-1}}=0$ or $1_{c\{u\}}1_{a_ib_i^{-1}}=1_{c\{u\}}$, $1_{c\{u\}}1_{e_jA_j}=0$ or $1_{c\{u\}}1_{e_jA_j}=1_{c\{u\}}$, and $1_{c\{u\}}1_{B_k}=0$ or $1_{c\{u\}}1_{B_k}=1_{c\{u\}}$. 
Then $$0\neq 1_{c\{u\}}x_0=\left(\sum\limits_{i:1_{c\{u\}}1_{a_ib_i^{-1}}\neq 0}\alpha_i+\sum\limits_{j:1_{c\{u\}}1_{e_jA_j}\neq 0}\beta_j+\sum\limits_{k:1_{c\{u\}}1_{B_k}\neq 0}\gamma_k\right)1_{c\{u\}}=h1_{c\{u\}}.$$
Therefore $(h1_{c\{u\}})\delta_0\in I$. Since $I$ is an ideal then $1_{c^{-1}}\delta_{c^{-1}}h1_{c\{u\}}\delta_0 1_c\delta_c =h\beta_{c^{-1}}(1_c1_{c\{u\}})\delta_0=(h1_u)\delta_0$ belongs to $I$.

Now assume that $1_dx_0\neq 0$ for some $d\in L$. Since $|d|>|a_i|$ then $1_d1_{a_ib_i^{-1}}=1_d$ or $1_d1_{a_ib_i^{-1}}=0$, for each $i\in \{1,...,m\}$, and similarly $1_d1_{e_j}A_j=1_d$ or $1_d1_{e_j}A_j=0$ for each $j\in \{1,...,n\}$. Moreover $1_d1_{B_k}=1_d$ if $s(d)\in B_k$, and $1_d1_{B_k}=0$ if $s(b)\notin B_k$, for each $k\in \{1,...,p\}$. Then we get that
$$0\neq 1_dx_0=\left(\sum\limits_{i:1_d1_{a_ib_i^{-1}}\neq 0}\alpha_i+\sum\limits_{j:1_d1_{e_jA_j}\neq 0}\beta_j+\sum\limits_{k:1_d1_{B_k}\neq 0}\gamma_k\right)1_d=h1_d,$$ and so $(h1_d)\delta_0\in I$. Hence $(h1_{r(d)})\delta_0=1_{d^{-1}}\delta_{d^{-1}}h1_d\delta_01_d\delta_d$  belongs to $I$. Then, for each vertex $w\in r(d)$, we get that $(h1_w)\delta_0=1_w\delta_0(h1_{r(d)})\delta_0$ belongs to $I$.
\fim

As a consequence of the above proposition we can provide a new proof of the Cuntz-Krieger Uniqueness Theorem for Leavitt   path algebras of ultragraphs.

\begin{corolario}\label{ck} Let $\GG$ be an ultragraph that satisfies Condition~$(L)$, let $R$ be commutative ring with a unit, and let $\pi:L_R(\GG)\rightarrow S$ be a homomorphism such that $\pi(r p_A)\neq 0$ for each $A\in \GG^0$ and non-zero $r\in R$. Then $\pi$ is injective. 
\end{corolario}
\demo 

Let $I=ker(\pi)$ and suppose that $I\neq 0$. Since $\GG$ satisfies Condition~$(L)$ then, by Theorem~\ref{maximalcommutative}, $D\delta_0$ is maximal commutative. Therefore, by \cite[Theorem 2.1]{johandanieldanilo}, $I\cap D\delta_0\neq 0$. Let $0\neq x_0\delta_0\in I\cap D_0\delta_0$. By Proposition~\ref{1vtheorem} there exist a non-zero $h\in R$, and an vertex $v$, such that $(h1_v)\delta_0\in I$, a contradiction. Therefore $ker(\pi)=0$.
\fim

As in the C* setting, the characterization of simplicity of ultragraph Leavitt path algebras rely on the notion of  hereditary and saturated collections. For the reader's convenience we recall these below. 

\begin{definicao} Let $\mathcal{G}$ be an ultragraph. A subcollection $H\subseteq\mathcal{G}^0$ is called hereditary if:
\begin{enumerate}
\item $s(e)\in H$ implies $r(e)\in H$, for each $e\in \mathcal{G}^1$;
\item $A\cup B\in H$, for all $A,B\in H$;
\item $A\in H, B\in\mathcal{G}^0$ and $B\subseteq A$ imply $B\in H$.
\end{enumerate}

Moreover, $H$ is called saturated if for any $v\in G^0$ with $0<|s^{-1}(v)|<\infty$, it holds that
$$\{r(e):e\in \mathcal{G}^1 \text{ and } s(e)=v\}\subseteq H \text{ implies } v\in H.$$
\end{definicao}

The next Lemma is key in the characterization of $\F$-simplicity in terms of existence of hereditary and saturated subcollections of $\GG^0$.

\begin{lema}\label{lemmaFinvariant} Let $R$ be a unital commutative domain and let $I$ be an $\F$-invariant ideal of $D_0$. Then the collection $$H=\{A\in \mathcal{G}^0: h1_A\in I \text{ for some non-zero } h\in R\}$$ is hereditary and saturated.
\end{lema}

\demo First we show that $H$ is hereditary. Let $e\in \mathcal{G}^1$ be such that $s(e)\in H$, and let $h\in R$ be a non-zero element such that $h1_{s(e)}\in I$. Then $h1_e=h1_e1_{s(e)}\in I\cap D_e$ and, since $I$ is $\F$-invariant, we have that $h1_{r(e)}=h\beta_{e^{-1}}(1_e)\in I$, and so $r(e)\in H$. Let $A,B\in H$, and let $h,k$ be non-zero elements in $R$ such that $h1_A\in I$ and $k1_B\in I$. Then $hk\neq 0$ since $R$ is a domain. Moreover, $hk1_{A\cup B}=hk1_A+hk1_B-hk1_A1_B\in I$ since $I$ is an ideal. Finally, let $A\in H$, and $B\in \mathcal{G}^0$ with $B\subseteq A$. Take a non-zero element $h\in R$ such that $h1_A\in I$. Note that $h1_B=h1_B1_A\in I$. Hence $B\in H$ and $H$ is hereditary.

Now we show that $H$ is saturated. Let $v\in G^0$ be such that $0<|s^{-1}(v)|<\infty$. Suppose that for each $e\in s^{-1}(v)$, it holds that $r(e)\in H$. Then for each $e\in s^{-1}(v)$ there is a non-zero $h_e\in R$ such that $h_e1_{r(e)}\in I$. Since $I$ is $\F$-invariant then $h_e1_e=h_e\beta_e(1_{e^{-1}})=\beta_e(h_e1_{r(e)})\in I$. Define $h=\prod\limits_{e\in s^{-1}(v)}h_e$, which is non-zero since $R$ is a domain. Then $h1_e\in I$ for each $e\in s^{-1}(v)$ and so $h1_v=\sum\limits_{e\in s^{-1}(v)}h1_e\in I$, from where we get that $v\in H$ and H is saturated. 

\fim

We can now describe the relation between $\F$-simplicity of $D$ and hereditary and saturated subcollections of $\mathcal{G}^0$.

\begin{teorema}\label{Fsimple} Let $R$ be a field. Then, the algebra $D$ is $\F$-simple if, and only if, the only hereditary and saturated subcollections of $\mathcal{G}^0$ are $\emptyset$ and $\mathcal{G}^0$.
\end{teorema}
\demo 

Suppose first that the only saturated and hereditary subcollections of $\mathcal{G}^0$ are $\emptyset$ and $\mathcal{G}^0$. Let $I\subseteq D$ be a non-zero, $\F$-invariant ideal. We show that $I=D$. Let $J$ be the set of all finite sums $\sum a_t\delta_t$, with $a_t\in D_t\cap I$. Notice that $J$ is is non-zero and is an ideal of $D\rtimes_\beta\F$, since $I$ is $F$-invariant. Then, by Proposition~\ref{1vtheorem}, there exists a $v\in G^0$, and a non-zero $h\in R$ such that $h1_v \delta_0 \in J$. Since $J\cap D_0\delta_0=I\delta_0$ then $h1_v\in I$. Let $H=\{A\in \mathcal{G}^0:h1_A\in I \text{ for some non-zero } h\in R\}$. By Lemma~\ref{lemmaFinvariant} $H$ is hereditary and saturated (and $H\neq \emptyset$ since $v\in H$), and hence $H=\mathcal{G}^0$. Then, for each $A\in \mathcal{G}^0$, there exists a non-zero element $h\in R$ such that $h1_A\in I$ and. Since $R$ is a field we have that $1_A\in I$, and it follows that $I=D_0$.

Now suppose that $D_0$ is $\F$-simple. Let $H\subseteq \mathcal{G}^0$ be nonempty, hereditary and saturated. We need to show that $H=\GG^0$. 

Let $I$ be the ideal in $D\rtimes_\beta\F$ generated by the set $\{1_A\delta_0:A\in H\}$, that is, $I$ is the linear span of all the elements of the form $a_r\delta_r1_A\delta_0a_s\delta_s$, with $r,s\in \F$, $a_r\in D_r$ and $a_s\in D_s$. Let $J=\{a:a\delta_0\in D\delta_0\cap I\}$, which is a non-zero ideal of $D$. Moreover, $J$ is $\F$ invariant, since if $a_t\in J\cap D_t$ then $a_t\delta_0\in I$ and $\beta_{t^{-1}}(a_t)\delta_0=1_{t^{-1}}\delta_{t^{-1}}a_t\delta_0 1_t\delta_t\in I$. Since $D$ is $\F$-simple then $J=D$.

Our next step is to show that $\{u\}\in H$, for each vertex $u\in G^0$.

Let $u\in G^0$. Then we can write $$1_u\delta_0=\sum\limits_tx_t\delta_t1_{A_t}\delta_0y_{t^{-1}}\delta_{t^{-1}}=\sum\limits_t \beta_t(\beta_{t^{-1}}(x_t)1_{A_t}y_{t^{-1}})\delta_0,$$ with $A_t\in H$. Multiplying the above equation by $1_{u}\delta_0$ we get that

\begin{equation}\label{eq.u}1_u=\sum\limits_{t\in T}1_u\beta_t(\beta_{t^{-1}}(x_t)1_{A_t}y_{t^{-1}}),
\end{equation}
where $T=\{t: 1_u\beta_t(\beta_{t^{-1}}(x_t)1_{A_t}y_{t^{-1}})\neq 0\}$.
In particular, for each $t\in T$ we have that $1_u1_t\neq 0$ and $1_{A_t}1_{t^{-1}}\neq 0$. 

If $u\in r(b)$, for some $b\in W$ with $\{s(b)\}\in H$, then $\{u\}\in H$ since $H$ is hereditary. If $0<|s^{-1}(u)|<\infty$, and $r(e)\in H$ for each $e\in s^{-1}(u)$, then $\{u\}\in H$ since $H$ is saturated. So we are left with the cases when there is no path $b$ with $\{s(b)\}\in H$ and $u\in r(b)$ and either $s^{-1}(u)=\emptyset$, $|s^{-1}(u)|=\infty$, or $0<|s^{-1}(u)|<\infty$ but $r(e)\notin H$ for some $e\in s^{-1}(u)$. Since there is no path $b\in W$ such that $\{s(b)\}\in H$ and $u\in r(b)$ then, for each $b\in W$, we get that $$1_u\beta_{b^{-1}}(\beta_b(x_{b^{-1}})1_Ay_b)=0$$
(notice that if $b\in W$ is such that $u\in r(b)$ then, since $H$ is hereditary, $s(b)\notin A$ and hence $1_A 1_{s(b)} y_b = 0$). So each non zero element $t\in T$ is of the form $t=ab^{-1}$, with $a\in W$ and $b\in W\cup \{0\}$.

{\it Case 1: $s^{-1}(u)=\emptyset$, and there is no path $b$ with $s(b)\in H$ and $u\in r(b)$}.

For each $t=ab^{-1}\in \F$ with $a\in W$ and $b\in W\cup\{0\}$, we have that $1_u1_t=0$, since $u$ is a sink, and so $t=ab^{-1}\notin T$. So $T=\{0\}$ and then $1_u=1_ux_01_{A_0}y_0$, with $A_0\in H$. Therefore $u\in A_0$ and so $\{u\}\in H$.

{\it Case 2: $|s^{-1}(u)|=\infty$, and there is no path $b$ with $\{s(b)\}\in H$ and $u\in r(b)$}.

 Suppose that $0\notin T$. Then each $t\in T$ is of the form $t=ab^{-1}$, with $a\in W$ and $b\in W\cup \{0\}$. Since $|s^{-1}(u)|=\infty$ then there exists $\xi\in X$ such that $s(\xi)\neq s(a)$ for each $ab^{-1}\in T$. So we get that $1=1_u(\xi)=\sum\limits_{t\in T}1_u\beta_t(\beta_{t^{-1}}(x_t)1_{A_t}y_{t^{-1}})(\xi)=0$, a contradiction. Hence $0\in T$, and so $1_ux_01_{A_0}y_0\neq 0$. Therefore $\{u\}\subseteq A_0\in H$ and, since $H$ is hereditary, we have that $\{u\}\in H$. 

Note that it follows from Case 1, Case 2, and by the fact that $H$ is hereditary, that if $u$ is a vertex such that $|s^{-1}(u)|=0$ or $|s^{-1}(u)|=\infty$ then $\{u\}\in H$.

{\it Case 3:  $0<|s^{-1}(u)|<\infty$, there is an edge $e\in s^{-1}(u)$ with $r(e)\notin H$, and there is no path $b$ with $\{s(b)\}\in H$ and $u\in r(b)$.}

Let us first prove the following claim:

{\it Claim: If $e$ is an edge such that $r(e)\notin H$ then there is a vertex $v\in r(e)$ such that $\{v\}\notin H$.}

Let $w=s(e)$. Notice that $\{w\}\notin H$, since $H$ is hereditary. Also note that there is no path $d$ with $s(d)\in \{H\}$ and $w\in r(d)$. Therefore, since $J=D$, proceeding as we did for $u$, we have that $$1_w=\sum\limits_{t\in S} 1_w\beta_t(\beta_{t^{-1}}(x'_t)1_{A_t}y'_{t^{-1}}),$$ where $A_t\in H$ for each $t\in S$, each 
$1_w\beta_t(\beta_{t^{-1}}(x'_t)1_{A_t}y'_{t^{-1}})$ is non zero, $0\notin S$ because $\{w\} \notin H$, and each $t$ is of the form $t=ab^{-1}$, with $a\in W$ and $b\in W\cup \{0\}$. 

For each $t=ab^{-1}\in T$ let $c_t=1_w\beta_t(\beta_{t^{-1}}(x'_t)1_{A_t}y'_{t^{-1}})$, so that 

\begin{equation}\label{eqw}1_w=\sum\limits_{t\in S}c_t.
\end{equation}

Since $1_w1_t\neq 0$ then $w=s(a)$ and, since $1_{A_t}1_{t^{-1}}\neq 0$, we have that $\{s(b)\}\subseteq A_t\in H$. Since $H$ is hereditary then $\{s(b)\}\in H$, and therefore $r(b)\in H$ and also $r(b)\cap r(a)\in H$. For $t=a\in W$ we get $A_t\cap r(t)\in H$.

By multiplying Equation (\ref{eqw}) on the left side by $1_{e^{-1}}\delta_{e^{-1}}$ and by $1_e\delta_e$ on the right side we get 

\begin{equation}\label{eqr(e)}
1_{r(e)}=\sum\limits_S \beta_{e^{-1}}(1_ec_t).
\end{equation}

Notice that for $t=a_1...a_{|a|}b^{-1}\in S$ with $a_1\neq e$ it holds that $\beta_{e^{-1}}(1_ec_t)=0$. Let $M=max\{|a|:ab^{-1}\in S \text{ and } a_1=e\}$, and let $S_i=\{ab^{-1}\in S:|a|=i \text{ and } a_1=e\}$, for $1\leq i\leq M$. In particular note that each element of $S_1$ is of the form $t=eb^{-1}$ with $b\in W\cup\{0\}$. 

If $e\notin S_1$ define $$A_1=\bigcup\limits_{ab^{-1}\in S_1}r(e)\cap r(b),$$ and if $e\in S_1$ define $$A_1=\left(\bigcup\limits_{ab^{-1}\in S_1, b\neq 0}r(e)\cap r(b)\right)\cup \left(r(e)\cap A_e\right).$$ Notice that $A_1\subseteq r(e)$ and that $A_1\in H$, since $r(e)\cap r(b)\in H$ for each $eb^{-1}\in S_1$ and $r(e)\cap A_e\in H$. 

From Equation (\ref{eqr(e)}) we get  $1_{r(e)}=\sum\limits_{i=1}^M\sum\limits_{t\in S_i}\beta_{e^{-1}}(1_ec_t)$. 

Now we show that $M>1$. Seeking a contradiction suppose $M=1$. Then we have that
$$1_{r(e)}=\sum\limits_{eb^{-1}\in S_1}\beta_{e^{-1}}(1_ec_{ab^{-1}}).$$

Since $A_1\subseteq r(e)$, $A_1\in H$, and $r(e)\notin H$, we have that $A_1$ is a proper subset of $r(e)$. So there is a vertex $v$ such that $v\in r(e)\setminus A_1$. Let $\xi \in X$ be such that $s(\xi)=v$ (notice that by the paragraph just above the statement of Case~3 $v$ is not a sink). Then for each $t=eb^{-1}\in S_1$ we get 
$$1_{eb^{-1}}(e\xi)=1_{e(r(e)\cap r(b))}(e\xi)=0,$$ since $r(e)\cap r(b)\subseteq A_1$.
Therefore
$$1=1_{r(e)}(\xi)=\sum\limits_{eb^{-1}\in S_1}\beta_{e^{-1}}(1_ec_{eb^{-1}})(\xi)=\sum\limits_{eb^{-1}\in S_1}1_{eb^{-1}}c_{eb^{-1}}(e\xi)=0.$$ a contradiction. Therefore $M>1$. 

Recall now that for each $ab^{-1}\in S_2\cup ...\cup S_M$ the element $a$ is of the form $a=a_1a_2...a_{|a|}=ea_2...a_{|a|}$. We want to show that $\{s(a_2)\}\notin H$ for some $ab^{-1}\in S_2\cup S_M$. Again seeking a contradiction, suppose that $\{s(a_2)\}\in H$, for each $ab^{-1}\in S_2\cup...\cup S_M$. Let $A_2$ be the set of all those vertices (the vertices $s(a_2)$). Notice that $A_2\in H$ (since we are supposing that each $\{s(a_2)\}\in H$ and $H$ is hereditary), and that $A_2\subseteq r(e)$ (since $s(a_2)\in r(a_1)=r(e)$). So we get that $A_1\cup A_2 \subseteq r(e)$ and, since $A_1\cup A_2 \in H$ and $r(e)\notin H$, there exist a vertex $v_0\in r(e)\setminus (A_1\cup A_2)$. Let $\xi \in X$with $s(\xi)=v_0$.

For each $eb^{-1}\in S_1$ we get $1_{eb^{-1}}(e\xi)=0$, since $s(\xi)\notin A_1$, and for each $ea_2...a_{|a|}b^{-1}\in S_2\cup ...\cup S_M$ we get $1_{ea_2...a_{|a|}b^{-1}}(e\xi)=0$, since $s(\xi)\neq s(a_2)$ (because $s(\xi)\notin A_2$). Therefore 
$$1=1_{r(e)}(\xi)=\sum\limits_{i=1}^M\sum\limits_{ab^{-1}\in S_i}\beta_{e^{-1}}(1_ec_{ab^{-1}})(\xi)=0,$$ a contradiction.

So there is an element $ab^{-1}\in S_2\cup...\cup S_M$ (where $a=ea_2...a_{|a|}$) with $\{s(a_2)\}\notin H$. Since $s(a_2)\in r(e)$, we proved the claim. 

Now we prove Case 3.

Firs write $1_u$ as in Equation (\ref{eq.u}), that is, 

$$1_u=\sum\limits_{t\in T}1_u\beta_t(\beta_{t^{-1}}(x_t)1_{A_t}y_{t^{-1}}),$$ where $1_u\beta_t(\beta_{t^{-1}}(x_t)1_{A_t}y_{t^{-1}})\neq 0$.

To show that $\{u\}\in H$ it is enough to show that $0\in T$, because in this case $0\neq 1_u1_{A_0}$, what implies that $u\in A_0$ and, since $A_0\in H$, then $\{u\}\in H$.

Suppose, by contradiction, that $0\notin T$. Then each $t\in T$ is of the form $t=ab^{-1}$ with $a\in W$ and $b\in W\cup \{0\}$. Recall that for each $t=ab^{-1}$ it holds that $r(a)\cap r(b)\in H$, and for $t=a$ it holds that $r(a)\cap A_a\in H$. 

Let $M=max\{|a|:ab^{-1}\in T, a\in W, b\in W\cup\{0\}\}$.

By hypothesis there is an edge $e_0\in s^{-1}(u)$ such that $r(e_0)\notin H$. By the previous claim, there is an vertex $v_1\in r(e_0)$ such that $\{v_1\}\notin H$. It follows from the paragraph right after Case 2 that $0<|s^{-1}(v_1)|<\infty$. Since $H$ is saturated there is an edge $e_1\in s^{-1}(v_1)$ such that $\{r(e_1)\}\notin H$. By applying the previous argument repeatedly we get a path $e_0...e_M$ such that $s(e_i)=v_i$, and $\{v_i\}\notin  H$, for each $i\in \{1,...,M\}$. Let $\xi \in X$ be such that $s(\xi)\in r(e_M)$. Then $e_0e_1...e_M\xi\in X$ and, for each $t=ab^{-1}\in T$, we get
$$1_{ab^{-1}}(e_0e_1...e_M\xi)=1_{a(r(a)\cap r(b))}(e_0e_1...e_M\xi)=0,$$ since $s(e_{|a|})\notin H$ and $r(a)\cap r(b)\in H$. The same holds for $t=a\in T$. So $1_t(e_0...e_M\xi)=0$ for each $t\in T$. Finally, we get that
$$1=1_u(e_0...e_M\xi)=\sum\limits_{t\in T}1_u\beta_t(\beta_{t^{-1}}(x_t)1_{A_t}y_{t^{-1}})(e_0...e_M\xi)=$$ $$=\sum\limits_{t\in T}1_u\beta_t(\beta_{t^{-1}}(x_t)1_{A_t}y_{t^{-1}})(e_0...e_M\xi)1_t(e_0...e_M\xi)=0,$$
a contradiction. Therefore $0\in T$ and Case~3 is proved.

So, we get that $\{u\}\in H$ for each $u\in \GG^0$.

To end the proof notice that, by \cite[Lemma~2.12]{Tom3}, any $A\in \GG^0$ can be written as $$\displaystyle
\bigcap_{e
\in X_1} r(e)
\cup \ldots 
\cup \bigcap_{e \in X_n} r(e) \cup F, $$where $X_1,
\ldots, X_n$ are finite subsets of $\mathcal{G}^1$, and $F$
is a finite subset of $G^0$. Since $H$ is hereditary and $\{s(e)\}\in H$, we have that $r(e)\in H$ for each $e\in \GG^1$. The result now follows from the fact that $H$ is hereditary.


\fim

We can now prove the simplicity criteria for the Leavitt path algebra of an ultragraph $\mathcal{G}$, $L_R(\mathcal{G})$, via partial skew group ring theory. 

\begin{teorema}\label{simlicitydescribed} Let $\mathcal{G}$ be an ultragraph and $R$ be a field. Then $L_R(\mathcal{G})$ is simple if, and only if, $\mathcal{G}$ satisfies condition $(L)$ and the unique saturated and hereditary subcollections of $\mathcal{G}^0$ are $\emptyset$ and $\mathcal{G}^0$. 
\end{teorema}
\demo 

By Theorem \ref{isom}, $L_R(\mathcal{G})$ and $D\rtimes_\beta\F$ are isomorphic algebras. By \cite[2.3]{johandanieldanilo}, the algebra $D\rtimes_\beta\F$ is simple if, and only if, $D$ is $\F$-simple and $D\delta_0$ is maximal commutative in $D\rtimes_\beta\F$. The result now follows from Theorems \ref{maximalcommutative} and \ref{Fsimple}.

\fim

In \cite[Teorem 3.11]{TomSimple} Tomforde gives a complete combinatorial description of ulgragraphs such that the associated ultragraph C*-algebra is simple. Since this description is obtained based only on the description of simplicity via hereditary and saturated collections the theorem above implies that we have the same description for $L_R(\GG)$. For reader's convenience we state the theorem below, but for this we need to recall a few definitions. 

For an ultragraph $\GG$, and $v,w\in G^0$, the notation $w\geq v$ means that there is a path $\alpha$ with $s(\alpha)=w$ and $v\in r(\alpha)$. Also, $G^0\geq \{v\}$ means that $w\geq v$ for each $w\in G^0$. The ultragraph $\GG$ is said to be cofinal if for each infinite path $\alpha=e_1e_2...$, and each vertex $v\in G^0$, there is an $i\in \N$ such that $v\geq s(e_i)$. Moreover, for $v\in G^0$ and $A\subseteq G^0$ we write $v\rightarrow A$ to mean that there are paths $\alpha_1,...,\alpha_n$ such that $s(\alpha_i)=v$, for all $1\leq i\leq n$, and $A\subseteq \bigcup\limits_{i=1}^nr(\alpha_i)$.

\begin{teorema} Let $\GG$ be an ultragraph and $R$ be a field. Then $L_R(\GG)$ is simple if and only if:

\begin{enumerate}
\item $\GG$ satisfies condition $(L)$ 
\item $\GG$ is cofinal
\item $G^0\geq \{v\}$ for every singular vertex $v\in G^0$
\item If $e\in \GG^1$ is an edge for which the set $r(e)$ is infinite, then for every $w\in G^0$ there exists a set $A_w\subseteq r(e)$ for which $r(e)\setminus A_w$ is finite and $v\rightarrow A_w$.
\end{enumerate}
\end{teorema}
\demo

The proof of this theorem relies only on the fact that the only hereditary and saturaded subcollections of $\GG^0$ are $\emptyset$ and $\GG^0$. So the proof given in \cite[Theorem~3.11]{TomSimple} applies. 

\fim

\section{Chain conditions}

In \cite{OinertChain} chain conditions are described for partial skew groupoid rings. As an application a new proof of the criteria for a Leavitt path algebra to be artinian is given. Namely, a Leavitt path algebra associated to a graph $E$ is artinian iff $E$ is finite and acyclic (A graph (ultragraph) is called acyclic if there are no cycles in the graph (ultragraph)). Building from the ideas in \cite{OinertChain} we show that this same criteria is true for ultragraph Leavitt path algebras. In our proof we will use that any ultragraph Leavitt path algebra of a finite acyclic ultragraph is isomorphic to a Leavitt path algebra of a finite acyclic graph, a result we state precisely below. 

Let $\mathcal{G}=(G^0, \mathcal{G}^1, r,s)$ be a finite ultragraph. Enumerate $G^0$, say $$G^0= \{v_1, \ldots, v_n\}.$$ Define a map $c:\GG^1 \rightarrow \{0,1\}^n$ by $c(e)=(y_i)$, where $y_i = \begin{cases} 1 &\mbox{if } v_i \in r(e)\\ 
0& \mbox{if } v_i \notin r(e) .\end{cases}$
Consider the graph $\mathcal{F}=(G^0, \mathcal{F}^1, r,s)$, where the set of edges $\mathcal{F}^1$ consists of all edges defined as follows: For each edge $e \in \GG^1$ and $i \in \{1,\ldots, n\}$ such that $c(e)_i =1$, let $f_{e_i}$ be the edge such that $s(f_{e_i})= s(e)$ and $r(f_{e_i}) = v_i$. We can now state the following proposition, a proof of which is left to the reader.

\begin{proposicao}\label{finite} Let $\GG$ be a finite ultragraph, that is, suppose that $G^0$ and $\GG^1$ are finite, and let $\mathcal{F}$ be the associated graph as defined above. Then $L_R(\GG)$ is isomorphic to $L_R(\mathcal{F})$. Furthermore, if $\GG$ is acyclic then $\mathcal{F}$ is acyclic. 
\end{proposicao}



We end the paper with the characterization of artinian ultragraph Leavitt path algebras. Recall that a ring is left (right) artinian if it satisfies the descending chain condition on left (right) ideals, and artinian if it is both left and right artinian.

\begin{teorema}\label{LPAthm}
Let $R$ be a field and let $\GG$ be an ultragraph.
Consider $L_R(\GG)$, the ultragraph Leavitt path algebra of $\GG$. Then the following five assertions are equivalent:
\begin{enumerate}[{\rm (i)}]
	\item\label{LPAthm:Econd} $\GG$ is finite and acyclic;
	\item\label{LPAthm:leftLPAartin} $L_R(\GG)$ is left artinian;
	\item\label{LPAthm:rightLPAartin} $L_R(\GG)$ is right artinian;
	\item\label{LPAthm:LPAartin} $L_R(\GG)$ is artinian;
	\item\label{LPAthm:LPAuss} $L_R(\GG)$ is unital and semisimple.
\end{enumerate}
\end{teorema}

\demo

All we need to prove is that \eqref{LPAthm:leftLPAartin}$\Rightarrow$\eqref{LPAthm:Econd}. The other implications follow from Proposition~\ref{finite} and \cite[Theorem~5.2]{OinertChain}.

\eqref{LPAthm:leftLPAartin}$\Rightarrow$\eqref{LPAthm:Econd}: The proof of this implication will follow closely the proof of \cite[Theorem~5.2]{OinertChain} for Leavitt path algebras. We include it here for completeness.

Suppose that $L_R(\GG) \cong D \rtimes_{\beta} \F$ is left artinian.
By \cite[Theorem~1.3]{OinertChain}, we get that $D_g=\{0\}$ for all but finitely many $g\in \F$, and $D$ is left artinian.

Assume that there exists an infinite path $p=e_1 e_2 e_3 \ldots$ in $\GG$. Then the ideals $D_{e_1}$, $D_{e_1e_2}$, $D_{e_1e_2 e_3}$, \ldots are all non-zero, a contradiction. Therefore there is no infinite path in $\GG$, and hence $\GG$ must be acyclic. 

Next we prove that $\GG$ is finite. Notice that if $G^0=\{v_1,v_2,v_3,\ldots\}$ is infinite then 
\begin{displaymath}
\oplus_{v\in E^0 \setminus \{v_1\}} L_R(\GG)v \supseteq \oplus_{v\in E^0 \setminus \{v_1,v_2\}}L_R(\GG)v \supseteq \oplus_{v\in E^0 \setminus \{v_1,v_2,v_3\}} L_R(\GG)v \supseteq \ldots
\end{displaymath}
is a descending chain of left ideals of $L_R(\GG)$ that never stabilizes (since every pair of vertices in $G^0$ are orthogonal idempotents).
Hence, $L_R(\GG)$ is not left artinian, a contradiction. Therefore $G^0$ is finite.

We finish the proof showing that $\GG^1$ is finite. Since $G^0$ is finite it is enough to prove that $G^0$ contains no infinite emitter.
Seeking a contradiction, suppose that there is a vertex $v\in G^0$ which is an infinite emitter.
Since $G^0$ is finite, there must exist some $u\in G^0$ such that
the set $I=\{e \in E^1 \mid s(e)=v \text{ and } u\in r(e) \}$ is infinite. If $u$ is a sink then $(u,u)\in X_{e^{-1}}$ for all $e\in I$, and hence $D_{e^{-1}}$ is non-zero for infinitely many $e\in I$, a contradiction. Suppose $u$ is not a sink. Then there exists a path $\eta\in X$ such that $s(\eta) = u$. Hence $X_{e^{-1}}$ contains $\eta$ for each $e\in I$. Therefore $D_{e^{-1}}$ is non-zero for infinitely many $e\in I$, a contradiction.

\fim

\vspace{1.5pc}

Daniel Gon\c{c}alves, Departamento de Matem\'{a}tica, Universidade Federal de Santa Catarina, Florian\'{o}polis, 88040-900, Brasil

Email: daemig@gmail.com

\vspace{0.5pc}
Danilo Royer, Departamento de Matem\'{a}tica, Universidade Federal de Santa Catarina, Florian\'{o}polis, 88040-900, Brasil

Email: daniloroyer@gmail.com
\vspace{0.5pc}

\end{document}